\RequirePackage{ifpdf}
\ifpdf 
\documentclass[pdftex]{sigma}
\else
\documentclass{sigma}
\fi

\begin{document}

\allowdisplaybreaks

\renewcommand{\PaperNumber}{064}

\FirstPageHeading

\ShortArticleName{Degenerate Series Representations of the $q$-Deformed
Algebra ${\rm so}'_q(r,s)$}

\ArticleName{Degenerate Series Representations\\ of the $\boldsymbol{q}$-Deformed
Algebra $\boldsymbol{{\rm so}'_q(r,s)}$}

\Author{Valentyna A. GROZA}
\AuthorNameForHeading{V.A. Groza}

\Address{National Aviation University, 1 Komarov Ave.,
03058 Kyiv, Ukraine}

\Email{\href{mailto:groza@i.com.ua}{groza@i.com.ua}}

\ArticleDates{Received January 26, 2007, in f\/inal form April
18, 2007; Published online May 02, 2007}

\Abstract{The $q$-deformed algebra ${\rm so}'_q(r,s)$  is a real
form of the $q$-deformed algebra $U'_q({\rm so}(n,\mathbb{C}))$,
$n=r+s$, which dif\/fers from the quantum algebra $U_q({\rm
so}(n,\mathbb{C}))$ of Drinfeld and Jimbo. We study representations
of the most degenerate series of the algebra ${\rm so}'_q(r,s)$. The
formulas of action of operators of these representations upon the
basis corresponding to restriction of representations onto the
subalgebra ${\rm so}'_q(r)\times {\rm so}'_q(s)$ are given. Most of
these representations are irreducible. Reducible representations
appear under some conditions for the parameters determining the
representations. All irreducible constituents which appear in
reducible representations of the degenerate series are found. All
$*$-representations of ${\rm so}'_q(r,s)$ are separated in the set
of irreducible representations obtained in the paper.}


\Classification{16B35; 16B37; 81R50}

\section{Introduction}

In this paper we consider most degenerate series representations of
the $q$-deformed algebra ${\rm so}'_q(r,s)$, which is a real form of
the complex $q$-deformed algebra $U'_q({\rm so} (n, \mathbb{C}))$
def\/ined in \cite{1}.  The algebra $U'_q({\rm so}(n,\mathbb{C}))$
dif\/fers from the quantum algebra $U_q({\rm so}(n,\mathbb{C}))$
def\/ined by Drinfeld~\cite{2} and Jimbo~\cite{3} (see also
\cite{KS}). Drinfeld and Jimbo def\/ined $U_q({\rm so}(n,\mathbb{C}))$
by means of Cartan subalgebras and root subspaces of the Lie algebra
${\rm so}(n,\mathbb{C})$. However, the Lie algebra ${\rm
so}(n,\mathbb{C})$ has a different structure based on the basis elements
$I_{k,k-1}=E_{k,k-1}-E_{k-1,k}$ (where $E_{is}$ is the matrix with
elements $(E_{is})_{rt}=\delta _{ir} \delta _{st})$. The
$q$-deformation of this structure leads to the algebra $U'_q({\rm
so}(n,\mathbb{C}))$, determined in \cite{1}. Later on it was shown
that this $q$-deformation of ${\rm so}(n,\mathbb{C})$ is very useful
in many directions of contemporary mathematics. Namely,
representations of the algebra $U'_q({\rm so}(n, \mathbb{C}))$ and
of its real forms are closely related to the theory of
$q$-orthogonal polynomials and $q$-special functions. Some ideas of
such applications can be found in \cite{4}.

The algebra $U'_q({\rm so} (n, \mathbb{C}))$ (especially its
particular case $U'_q({\rm so}(3,\mathbb{C}))$) is related to the
algebra of observables in 2+1 quantum gravity on the Riemmanian
surfaces (see, for example, \cite{5}). A~quantum analogue of the
Riemannian symmetric space $SU(n)/SO(n)$ is constructed by means of
the algebra $U'_q({\rm so} (n, \mathbb{C}))$~\cite{6}. It is clear
that a construction of a quantum analogue of some pseudo-Riemannian
symmetric spaces is connected with the $q$-deformed algebra ${\rm
so}'_q(r,s)$.

A $q$-analogue of the theory of harmonic polynomials ($q$-harmonic
polynomials on quantum vector space ${\mathbb{R}}_q^n$) is constructed by
using the algebra $U'_q({\rm so}
(n, \mathbb{C}))$. In particular, a $q$-analogue
of separations of variables for the $q$-Laplace operator on ${\mathbb{R}}_q^n$
is given by means of this algebra and its subalgebras (see \cite{7,8}).
The algebra $U'_q({\rm so} (n, \mathbb{C}))$ also appears in the theory of links
in the algebraic topology \cite{9}.

The representation theory of the $q$-deformed algebra
$U'_q({\rm so} (n, \mathbb{C}))$ dif\/fers from that for
the Drinfeld--Jimbo algebra $U_q({\rm so} (n, \mathbb{C}))$.
One of these dif\/ferences is related to the fact that
the Drinfeld--Jimbo algebra $U_q({\rm so} (n, \mathbb{C}))$
admits the inclusion
\[
U_q({\rm so}(n,{\mathbb{C}}))\supset U_q({\rm so}(n-2,{\mathbb{C}}))
\]
and does not admit the inclusion
\[
U_q({\rm so}(n,{\mathbb{C}}))\supset U_q({\rm so}(n-1,{\mathbb{C}})).
\]
The algebra $U'_q({\rm so} (n, \mathbb{C}))$ admits such an
inclusion. This allows to construct Gel'fand--Tsetlin bases for
f\/inite dimensional representations of $U'_q({\rm so} (n, \mathbb{C}))$
(see \cite{1}).

As in the case of real forms of Drinfeld--Jimbo quantum algebras
(see \cite{11a,12a,13a}) the real form
${\rm so}'_q(r,s)$ of $U'_q({\rm so}(r+s,\mathbb{C}))$ is def\/ined by
means of introducing a
$*$-operation into $U'_q({\rm so}(r+s,\mathbb{C}))$. When $q\to 1$ then
the $q$-deformed algebra ${\rm so}'_q(r,s)$ turns into the universal
enveloping algebra $U({\rm so}_{r,s})$ of the Lie algebra ${\rm so}_{r,s}$
which corresponds to the pseudo--orthogonal Lie group $SO_0(r,s)$.

Representations of the algebra ${\rm so}'_q(r,s)$, considered in this paper,
are given by one continuous parameter. These representations are
$q$-deformations of the representations of the classical Lie group $SO_0(r,s)$
considered in \cite{10,11} (see also \cite{12,13}). We derive several series of
$*$-representations of the algebra ${\rm so}'_q(r,s)$. As in the case of the
quantum algebra $U_q({\rm su}_{1,1})$, the algebra ${\rm so}'_q(r,s)$
has the so-called strange series of $*$-representations, which is absent in
the case of the Lie group $SO_0(r,s)$. When $q\to 1$ then
matrix elements of operators of the strange series representations tend to
the inf\/inity and representations become senseless.

Everywhere below we consider that $q$ is a
positive number. We also suppose that $r>2$ and $s>2$. Representations of degenerate
series of the algebra ${\rm so}'_q(r,1)$ were considered in \cite{14}.
Representations of the algebra ${\rm so}'_q(r,2)$ were studied in \cite{15}.
In fact, we generalize a part of the results of papers \cite{14,15}.

It is well-known that the algebra $U'_q({\rm so}(n,\mathbb{C}))$ has
f\/inite dimensional irreducible representations of two types:
representations of the classical type (at $q\to 1$ they tend to
the corresponding representations of the Lie algebra
${\rm so}(n,\mathbb{C})$) and representations of the non-classical
type (there exists no analogue of these representations in the
case of ${\rm so}(n,\mathbb{C})$). The algebra ${\rm so}'_q(r,s)$
has no degenerate series representations of the non-classical type.
The reason is that the ``compact'' algebras
${\rm so}'_q(n)$ have no degenerate irreducible representations
of the non-classical type.

\section[The $q$-deformed algebra ${\rm so}'_q(r,s)$]{The $\boldsymbol{q}$-deformed
algebra $\boldsymbol{{\rm so}'_q(r,s)}$}

 The algebra ${\rm so}'_q(r,s)$ is a real form of the
$q$-deformed algebra
$U'_q({\rm so}(r+s, \mathbb{C}))$ which is separated by the $*$-operation.
The algebra $U'_q({\rm so}(r+s, \mathbb{C}))$ is def\/ined in
\cite{1}.

The classical universal enveloping algebra $U({\rm so} (n,\mathbb{C}))$
is generated by the
elements $I_{i,i-1}$, $i=2,3, \ldots ,n$, that satisfy the relations
\begin{gather}\label{1}
I_{i,i-1}I^2_{i+1,i}-2I_{i+1,i}I_{i,i-1}I_{i+1,i}+I^2_{i+1,i}I_{i,i-1}
=-I_{i,i-1},
\\
\label{2}
I^2_{i,i-1}I_{i+1,i}-2I_{i,i-1}I_{i+1,i}I_{i,i-1}+I_{i+1,i}I^2_{i,i-1}
=-I_{i+1,i},
\\
\label{3}
[I_{i,i-1},I_{j,j-1}]=0,\qquad \vert i-j\vert >1
\end{gather}
(they follow from the well-known commutation relations for the generators
$I_{ij}$ of the Lie algebra ${\rm so} (n,\mathbb{C}))$.
In the approach to the $q$-deformed orthogonal algebra of paper \cite{1}, a
$q$-deformation of the associative algebra $U({\rm so} (n,\mathbb{C}))$
is def\/ined by deforming the relations \eqref{1}--\eqref{3}.
These $q$-deformed relations are of the form
\begin{gather}\label{4}
I_{i,i-1}I_{i+1,i}^2-aI_{i+1,i}I_{i,i-1}I_{i+1,i}
+I^2_{i+1,i}I_{i,i-1}=-I_{i,i-1},
\\
\label{5}
I^2_{i,i-1}I_{i+1,i}-aI_{i,i-1}I_{i+1,i}I_{i,i-1}
+I_{i+1,i}I^2_{i,i-1}=-I_{i+1,i},
\\
\label{6}
[I_{i,i-1},I_{j,j-1}]=0,\ \ \vert i-j\vert >1,
\end{gather}
where $a=q^{1/2}{+}q^{-1/2}$ and $[\cdot ,\cdot ]$ denotes
the usual commutator. Obviously, in the
limit $q {\to} 1$ formulas \eqref{4}--\eqref{6} give relations
\eqref{1}--\eqref{3}. Remark that
relations \eqref{4} and \eqref{5} dif\/fer from the $q$-deformed Serre relations in the
approach of Jimbo and Drinfeld to the quantum algebras
$U_q({\rm so} (n,\mathbb{C}))$ by
appearance of nonzero right hand sides and by possibility of reduction
\[
U'_q({\rm so}(n,\mathbb{C}))\supset U'_q({\rm so}(n-1,\mathbb{C})).
\]
Below, by the algebra $U'_q({\rm so} (n,\mathbb{C}))$ we mean the $q$-deformed
algebra def\/ined by formulas \eqref{4}--\eqref{6}.

The ``compact'' real form ${\rm so}'_q(n)$ of the algebra $U'_q({\rm so}
(n,\mathbb{C}))$ is def\/ined by the involution given~as
\begin{gather}\label{8}
I^*_{i,i-1}=-I_{i,i-1},\qquad i=2,3,\ldots ,n.
\end{gather}
The ``noncompact'' real form ${\rm so}'_q(r,n-r)$ of $U'_q({\rm so} (n,
\mathbb{C}))$ is determined by the involution
\begin{gather}\label{9}
I^*_{i,i-1}=-I_{i,i-1},\qquad i\ne r+1,\qquad I^*_{r+1,r}=I_{r+1,r}.
\end{gather}
It would be more correct to use the notation $U'_q({\rm so}_q(r,n-r))$
for ${\rm so}'_q(r,n-r)$. We use the last notation since it is simpler.

The $q$-deformed algebra ${\rm so}'_q(n)$
contains the subalgebra ${\rm so}'_q(n-1)$. This fact allows
us to consider Gel'fand--Tsetlin bases of carrier spaces of
representations of ${\rm so}'_q(n)$ \cite{1}.
The $q$-deformed algebra ${\rm so}'_q(n-r,r)$ contains the subalgebra
${\rm so}'_q(n-r)\times {\rm so}'_q(r)$.

\section[Representations of ${\rm so}'_q(n)$]{Representations of $\boldsymbol{{\rm so}'_q(n)}$}

Irreducible representations of the algebra ${\rm so}'_q(r,s)$ are
described by means of f\/inite dimensional irreducible representations of
the subalgebras ${\rm so}'_q(r)$ and ${\rm so}'_q(s)$. Therefore, we
describe representations of ${\rm so}'_q(n)$ which will be used below.

Irreducible f\/inite dimensional representations of the algebra
${\rm so}'_q(3)$ (belonging to the classical type)
are given by integral or half-integral nonnegative
number $l$. We denote these representations by $T_l$. The carrier space of
the representation $T_l$ has the orthonormal basis
\[
|m\rangle ,\qquad  m=-l,-l+1,\ldots ,l,
\]
and the operators $T_l(I_{21})$ and $T_l(I_{32})$ act upon this basis as (see \cite{14})
\begin{gather}\label{11}
T_l(I_{21})|m\rangle ={\rm i}[m]_q|m\rangle ,\ \ \ {\rm i}=\sqrt {-1},
\\
T(I_{32})|m\rangle =d(m)\left( [l-m]_q[l+m+1]_q\right)^{1/2}|m+1\rangle
\nonumber\\
\phantom{T(I_{32})|m\rangle =}{}
-d(m-1)\left( [l+m]_q[l-m+1]_q\right)^{1/2}|m-1\rangle ,\label{12}
\end{gather}
where
\[
d(m)=\left( {[m]_q[m+1]_q\over [2m]_q[2m+2]_q}\right) ^{1/2}
\]
and $[b]_q$ is a $q$-number def\/ined by the formula
\[
[b]_q={q^{b/2}-q^{-b/2}\over q^{1/2}-q^{-1/2}}.
\]

Let us describe f\/inite dimensional irreducible representations of
the algebra ${\rm so}'_q(n)$, $n>3$, which are of class 1 with respect to the
subalgebra ${\rm so}'_q(n-1)$~\cite{8}. As in the classical case, these
representations are given by highest weights $(m_n,0,\dots ,0)$, where
$m_n$ is a nonnegative integer. We denote these representations by $T_{m_n}$.
Under restriction to the subalgebra ${\rm so}'_q(n{-}1)$, the
representation $T_{m_n}$ contains (with unit multiplicity)
those and only those irreducible representations~$T_{m_{n-1}}$ of this subalgebra for which we have
\[
m_n\ge m_{n-1}\ge 0.
\]
Exactly in the same way as in the case of the classical group $SO(n)$, we
introduce the Gel'fand--Tsetlin basis of the carrier space of
the representation $T_{m_n}$ by using the successive reduction
\[
{\rm so}'_q(n)\supset {\rm so}'_q(n-1)\supset {\rm so}'_q(n-2)\supset
\cdots \supset {\rm so}'_q(3)\supset {\rm so}'_q(2).
\]
We denote the basis elements of this space by
\[
|m_n,m_{n-1},m_{n-2},\dots ,m_3,m_2\rangle ,
\]
where $m_n\ge m_{n-1}\ge m_{n-2}\ge \cdots \ge |m_2|$ and $m_{n-i}$ determines
the representation $T_{m_{n-i}}$ of ${\rm so}'_q(n-i)$. With respect to this basis
the operator $T(I_{n,n-1})$ of the representation $T_{m_n}$ of
${\rm so}'_q(n)$ is given by the formula
\begin{gather}
T_{m_n}(I_{n,n-1})|m_n,m_{n-1},\dots ,m_2\rangle\nonumber\\
\qquad
{}=([m_n+m_{n-1}+n-2]_q[m_n-m_{n-1}]_q)^{1/2}R(m_{n-1})|m_n,m_{n-1}+1,\dots ,m_2\rangle\label{13}\\
{}\qquad -([m_n+m_{n-1}+n-3]_q[m_n-m_{n-1}+1]_q)^{1/2}
R(m_{n-1}-1)
|m_n,m_{n-1}-1,\dots ,m_2\rangle ,\nonumber
\end{gather}
where
\[
R(m_{n-1})=\left( {[m_{n-1}+m_{n-2}+n-3]_q[m_{n-1}-m_{n-2}+1]_q\over
[2m_{n-1}+n-3]_q[2m_{n-1}+n-1]_q}\right) ^{1/2}.
\]
The other operators $T(I_{i,i-1})$ are given by the same formulas
with the corresponding change for $m_n$ and $m_{n-1}$ or at
$i=3,\, 2$ by the formulas for irreducible representations of the algebra
${\rm so}'_q(3)$, described above.

The representations $T_{m_n}$ are characterized by the property that under
restriction to the subalgebra ${\rm so}'_q(n-1)$ the restricted
representations contain a trivial irreducible representation of
${\rm so}'_q(n-1)$ (that is, a representation with highest weight
$(0,0,\dots,0)$). This is why one says that the representation $T_{m_n}$ is
of class 1 with respect to ${\rm so}'_q(n-1)$. The irreducible
representations~$T_{m_n}$ exhaust all irreducible
representations of class 1 of the algebra ${\rm so}'_q(n)$.

\section{Representations of the degenerate principal series}

We shall consider inf\/inite dimensional representations of the associative
algebra ${\rm so}'_q(r,s)$. Moreover, we admit representations by
unbounded operators. There exist non-equivalent def\/initions of
representations of associative algebras by unbounded or bounded
operators (see \cite{17,18}). In order to have a natural
def\/inition of a representation of ${\rm so}'_q(r,s)$ we take into account
the following items:
 \smallskip

(I) We shall deal also with $*$-representations (it is well-known
that these representations
are an analogue of unitary representations of Lie groups). Therefore,
for each representation operator there should exist an adjoint operator. This
means that a representation space have to be def\/ined on a Hilbert space.
 \smallskip

(II) Unbounded operators cannot be def\/ined on the whole Hilbert space.
However, existence of an adjoint operator $A^*$ to an unbounded operator
$A$ means that the operator $A$ must be def\/ined on an everywhere dense
subspace in the Hilbert space.
 \smallskip

(III) In order to be able to consider products of representation operators,
there must exist an everywhere dense subspace of the Hilbert space
which enter to a domain of def\/inition of each representation operator.
 \smallskip

Therefore, we give the following def\/inition of a representation of
${\rm so}'_q(r,s)$.
A representation $T$ of the associative algebra ${\rm so}'_q(r,s)$
is an algebraic homomorphism from ${\rm so}'_q(r,s)$
into an algebra of linear (bounded or unbounded) operators on a Hilbert
space $\mathcal{H}$ for which the following conditions are fulf\/illed:

\smallskip

a) the restriction of $T$ onto the ``compact'' subalgebra
${\rm so}'_q(r)\times {\rm so}'_q(s)$ decomposes into a direct sum of
its f\/inite dimensional irreducible representations (given by highest
weights) with f\/inite multiplicities;
\smallskip

(b) operators of a representation $T$ are determined on an everywhere
dense subspace $\mathcal{W}$ of $\mathcal{H}$,
containing all subspaces which are carrier
spaces of irreducible f\/inite dimensional representations of
${\rm so}'_q(r)\times {\rm so}'_q(s)$ from the restriction of $T$.
\smallskip

   In other words, our representations of ${\rm so}'_q(r,s)$
are Harish-Chandra modules of ${\rm so}'_q(r,s)$ with respect to
${\rm so}'_q(r)\times {\rm so}'_q(s)$.

There exist dif\/ferent non-equivalent def\/initions of irreducibility
of representations of associative algebras by unbounded operators
(see \cite{17}). Since unbounded representation operators are not
def\/ined on all elements of the Hilbert space, then we cannot def\/ine
irreducibility as in the f\/inite dimensional case. A natural
def\/inition is the following one.
A representation $T$ of ${\rm so}'_q(r,s)$ on $\mathcal{H}$ is called
irreducible if $\mathcal{H}$ has no non-trivial invariant subspaces such that
its closure does not coincide with $\mathcal{H}$. If operators of a
representation $T$ obey the relations
\begin{gather}\label{14}
T(I_{i,i-1})^*=-T(I_{i,i-1}),\qquad i\ne r+1,\qquad
T(I_{r+1,r})^*=T(I_{r+1,r})
\end{gather}
(compare with formulas \eqref{9}) on a common domain $\mathcal{W}$ of def\/inition, then
$T$ is called a $*$-representation.

 To def\/ine a representation $T$ of ${\rm so}'_q(r,s)$ it is
suf\/f\/icient to give the operators $T(I_{i,i-1})$, $i=2,3,\dots ,r+s$,
satisfying relations \eqref{4}--\eqref{6} on a common domain of
def\/inition. Let us def\/ine representations of
${\rm so}'_q(r,s)$ belonging to the {\it degenerate principal series}.
They are given by a~complex number $\lambda$ and a number $\epsilon \in
\{ 0, 1\}$. We denote the corresponding representation by
$T_{\epsilon \lambda }$. The space ${\mathcal{H}}(T_{\epsilon \lambda })$ of the
representation $T_{\epsilon \lambda }$ is an orthogonal sum of the subspaces
${\mathcal{V}}(m,{\bf 0}; m',{\bf 0})$, which are the carrier spaces of the f\/inite
dimensional representations of
${\rm so}'_q(r)\times {\rm so}'_q(s)$  with highest weights
$(m,{\bf 0}; m',{\bf 0})$ such that
\[
m+m'\equiv \epsilon \ ({\rm mod}\ 2).
\]
Here $(m,{\bf 0})$ and $(m',{\bf 0})$, $m\ge 0$, $m' \ge 0$, are highest weights of
irreducible representations of the subalgebras ${\rm so}'_q(r)$ and
${\rm so}'_q(r)$, respectively, and {\bf 0} denotes the set of zero
coordinates (in the case of the subalgebra ${\rm so}'_q(3)$, {\bf 0}
must be omitted). We assume that the basis vectors from \eqref{13} are orthonormal
in ${\mathcal{V}}(m,{\bf 0}; m',{\bf 0})$. Therefore, we have
\begin{gather}\label{15}
{\mathcal{H}}(T_{\epsilon \lambda })=\bigoplus _{m+m'\equiv \epsilon \; ({\rm mod}\ 2)}
{\mathcal{V}}(m,{\bf 0}; m',{\bf 0}),
\end{gather}
where we suppose that the sum means a closure of the corresponding linear
span. A linear span of the subspaces ${\mathcal{V}}(m,{\bf 0}; m',{\bf 0})$
determines an everywhere dense subspace on which all operators
of the representation $T_{\epsilon \lambda }$ are def\/ined.
Recall that we suppose that $r>2$ and $s>2$.

   We choose in the subspaces $\mathcal{V}(m,{\bf 0}; m',{\bf 0})$ the orthonormal
bases which are products of the bases corresponding to irreducible
representations of the subalgebras ${\rm so}'_q(r)$ and
${\rm so}'_q(s)$ introduced in Section 3 (Gel'fand--Tsetlin bases).
Elements of such bases are labeled by double Gel'fand--Tsetlin patterns
which will be denoted as
\begin{gather}\label{16}
|M\rangle =|m,k,j,\dots ;m',k',j',\dots \rangle .
\end{gather}
It is clear
that entries of patterns \eqref{16} obey the following conditions:
\begin{gather}\label{17}
m\ge k\ge j\ge \cdots ,\qquad m'\ge k'\ge j'\ge \cdots
\end{gather}
if $r>3$ and $s>3$. Thus, elements of the orthonormal basis of the carrier
space of the representation $T_{\epsilon \lambda }$ are labeled by all
patterns \eqref{16} satisfying betweenness conditions \eqref{17} and the equality
$m+m'\equiv \epsilon $ (mod 2).

   In order to def\/ine the representations $T_{\epsilon \lambda }$
we give explicit formulas for the operators
$T_{\epsilon \lambda }(I_{i,i-1}),\!$ $i=2,3,\dots ,r+s$. The operators
$T_{\epsilon \lambda }(I_{i,i-1})$, $i=2,3,\dots ,r$, act upon basis
elements \eqref{16} by the formulas of Section 3 as operators of the
corresponding irreducible representations of the subalgebra
${\rm so}'_q(r)$. It is clear that these operators act only upon
entries $k$, $j$, $\dots$ of vectors \eqref{16} and do not change entries
$m$, $m'$, $k'$, $j'$, $\dots$. The operators
$T_{\epsilon \lambda }(I_{i,i-1})$, $i=r+2, r+3,\dots ,r+s$, act
upon basis elements \eqref{16} by formulas of Section 3 as operators of corresponding
irreducible representations of the subalgebra ${\rm so}'_q(s)$.
The operator $T_{\epsilon \lambda }(I_{r+1,r})$ acts upon vectors \eqref{16} by the
formula
\begin{gather}
T_{\epsilon \lambda }(I_{r+1,r})|m,k,j,\dots ;m',k',j',\dots \rangle
\nonumber\\
\qquad{}=K_mL_{m'}[\lambda +m+m']_q|m+1,k,j,\dots ;m'+1,k',j',\dots \rangle
\nonumber\\
\qquad{}
-K_mL_{m'-1}[\lambda +m-m'-s+2]_q|m+1,k,j,\dots ;m'-1,k',j',\dots \rangle
\nonumber\\
  \qquad{} +K_{m-1}L_{m'}[\lambda -m+m'-r+2]_q|m-1,k,j,\dots ;m'+1,k',j',\dots \rangle
\nonumber\\
  \qquad{}-K_{m-1}L_{m'-1}[\lambda -m-m'-r-s+4]_q|m-1,k,j,\dots ;m'-1,k',j',\dots \rangle ,\label{18}
\end{gather}
where
\begin{gather}\label{19}
K_m=\left( {[m-k+1]_q[m+k+r-2]_q\over [2m+r]_q[2m+r-2]_q}\right) ^{1/2},
\\
\label{20}
L_{m'}=\left( {[m'-k'+1]_q[m'+k'+s-2]_q\over [2m'+s]_q[2m'+s-2]_q}\right) ^{1/2}.
\end{gather}
So, this operator changes only the entries $m$ and $m'$ in vectors \eqref{16}.

To prove that these operators really give a representation of the algebra
${\rm so}'_q(r,s)$ we have to verify that the relations
\begin{gather}
T_{\epsilon \lambda}(I_{r+1,r})^2T_{\epsilon \lambda}(I_{r,r-1})
-\big(q^{1\over 2}+q^{-{1\over 2}}\big)T_{\epsilon \lambda}(I_{r+1,r})T_{\epsilon \lambda}(I_{r,r-1})
T_{\epsilon \lambda}(I_{r+1,r})\nonumber\\
\qquad{}
+T_{\epsilon \lambda}(I_{r,r-1})T_{\epsilon \lambda}(I_{r+1,r})^2=
-T_{\epsilon \lambda}(I_{r,r-1}),\label{21}
\\
T_{\epsilon \lambda}(I_{r+1,r})T_{\epsilon \lambda}(I_{r,r-1})^2
-\big(q^{1\over 2}+q^{-{1\over 2}}\big)T_{\epsilon \lambda}(I_{r,r-1})T_{\epsilon \lambda}(I_{r+1,r})
T_{\epsilon \lambda}(I_{r,r-1})\nonumber\\
    \qquad{}+T_{\epsilon \lambda}(I_{r,r-1})^2T_{\epsilon \lambda}(I_{r+1,r})=
-T_{\epsilon \lambda}(I_{r+1,r}),\label{22}
\end{gather}
as well as relations \eqref{21} and \eqref{22}, in which $T_{\epsilon \lambda}(I_{r,r-1})$
is replaced by $T_{\epsilon \lambda}(I_{r+2,r+1})$, and the relations
\begin{gather}\label{23}
[T_{\epsilon \lambda}(I_{i,i-1}),T_{\epsilon \lambda}(I_{j,j-1})]=0,
\qquad |i-j|>1,
\end{gather}
where $[\cdot ,\cdot ]$ is the usual commutator, are fulf\/illed.

Fulf\/ilment of these relations can be shown by a direct calculation. Namely,
we act by both their  parts upon vector \eqref{16}, then collect coef\/f\/icients at
the same resulting basis vectors and show that the relations obtained are
correct. We do not give here these direct calculations.

Irreducibility of representations $T_{\epsilon \lambda}$
is studied by means of the following proposition.

\begin{proposition} The representation $T_{\epsilon \lambda}$
is irreducible if each of the numbers $[\lambda +m+m']_q$,
$[\lambda +m-m'-s+2]_q$, $[\lambda -m+m'-r+2]_q$,
$[\lambda -m-m'-r-s+4]_q$ in \eqref{18} vanishes only if
the vector on the right hand side of \eqref{18} with the coef\/f\/icient,
containing this number, does not
belong to the Hilbert space
${\mathcal{H}}(T_{\epsilon \lambda})$.
\end{proposition}

This proposition is proved in the same way as the corresponding proposition
in the classical case (see Proposition 8.5 in Section 8.5 of~\cite{12}). Thus, the main role under studying irreducibility of
the representations $T_{\epsilon \lambda}$ is played by the coef\/f\/icients
at vectors in \eqref{18} numerated in Proposition~1.

\section{Irreducibility}

To study the representations $T_{\epsilon \lambda }$ we take into account
properties of the function
\[
w(z)=[z]_q={q^{z/2}-q^{-z/2}\over q^{1/2}-q^{-1/2}}
={e^{hz/2}-e^{-hz/2}\over e^{h/2}-e^{-h/2}}.
\]
where $q=e^h$. Namely, we have
\[
w(z+4\pi k{\rm i}/h)=w(z),\qquad w(z+2\pi k{\rm i}/h)=-w(z),\qquad
k\ {\rm are\ odd\ integers}.
\]
These relations mean that the following proposition holds:

\begin{proposition} For arbitrary $\lambda$ the pairs of representations
$T_{\epsilon \lambda }$ and $T_{\epsilon ,\lambda +4\pi k{\rm i}/h}$ are coinciding
and the pairs of representations $T_{\epsilon \lambda }$ and
$T_{\epsilon ,\lambda +2\pi k{\rm i}/h}$ are equivalent.
\end{proposition}

Therefore, we may
consider only representations $T_{\epsilon \lambda }$ with
$0\le {\rm Im}\, \lambda <2\pi /h$.

   Most of the representations $T_{\epsilon \lambda }$ are irreducible.
Nevertheless, some of them are reducible. As in the classical case (see
\cite{10,11,12,13}),
reducibility appears because of vanishing of some of the coef\/f\/icients
$[\lambda +m+m']_q$, $[\lambda +m-m'-s+2]_q$, $[\lambda -m+m'-r+2]_q$,
$[\lambda -m-m'-r-s+4]_q$ in formula \eqref{18}. Using this fact we derive
(in the same way as in the classical case; see \cite{12}) the
following theorem:

\begin{theorem} If $r$ and $s$ are both even, then the representation
$T_{\epsilon \lambda }$ is irreducible if and only if $\lambda$ is not
an integer such that $\lambda \equiv \epsilon $ $({\rm mod}\; 2)$.  If one of the
numbers $r$ and $s$ is even and the other one is odd, then the representation
$T_{\epsilon \lambda }$ is irreducible if and only if $\lambda $ is not
an integer. If~$r$ and~$s$ are both odd, then the representation
$T_{\epsilon \lambda }$ is irreducible if and only if $\lambda$ is not an
integer or if $\lambda$ is an integer such that
\[
\lambda \equiv \epsilon \ ({\rm mod}\ 2),\qquad 0<\lambda <\tfrac 12 (r+s)-2.
\]
\end{theorem}

Irreducible representations $T_{\epsilon \lambda }$ admit additional equivalence
relations.

\begin{proposition} The pairs of irreducible representations $T_{\epsilon \lambda }$
and $T_{\epsilon ,r+s-\lambda -2}$ are equivalent.
\end{proposition}

 This equivalence is proved
by constructing an explicit form of intertwining operators. These operators
are diagonal in the basis \eqref{16} and can be evaluated exactly in the same way
as in the classical case.

   Sometimes it is useful to have the representations $T_{\epsilon \lambda }$
in somewhat dif\/ferent basis, namely, in the basis
\[
 |m,k,j,\dots ;\ m',k',j',\dots \rangle '
\]
which is related to the basis \eqref{16} by the formulas
\begin{gather*}
|m_0+\epsilon +i,k,j,\dots ; m_0-i,k',j',\dots \rangle
\\
\qquad{}
={\prod\limits_{t=1}^{m_0}[-\lambda +\epsilon +r+s+2t-4]_q^{1/2}\over
\prod\limits_{t=1}^{m_0}[\lambda +\epsilon +2t-2]_q^{1/2}}
{\prod\limits_{t=1}^i[-\lambda +\epsilon +r+2t-2]_q^{1/2}\over
\prod\limits_{t=1}^i[\lambda +\epsilon -s+2t]_q^{1/2}}
\\ \qquad{}
\times |m_0+\epsilon +i,k,j,\dots ; m_0-i,k',j',\dots \rangle ',
\\
|m_0+\epsilon -i,k,j,\dots ; m_0+i,k',j',\dots \rangle
\\   \qquad{}
={\prod\limits_{t=1}^{m_0}[-\lambda +\epsilon +r+s+2t-4]_q^{1/2}\over
\prod\limits_{t=1}^{m_0}[\lambda +\epsilon +2t-2]_q^{1/2}}
{\prod\limits_{t=1}^i[\lambda +\epsilon -s-2t+2]_q^{1/2}\over
\prod\limits_{t=1}^i[-\lambda +\epsilon +r-2t]_q^{1/2}}
\\
  \qquad{}
\times |m_0+\epsilon -i,k,j,\dots ; m_0+i,k',j',\dots \rangle '.
\end{gather*}
In the new basis the operators
$T_{\epsilon \lambda }(I_{i,i-1})$, $i\ne r+1$, have the same form as in the
basis \eqref{16}, and the operator $T_{\epsilon \lambda }(I_{r+1,r})$ is of the
form
\begin{gather}
T_{\epsilon \lambda }(I_{r+1,r})|m,k,j,\dots ;m',k',j',\dots \rangle '
\nonumber\\
\qquad{}
=K_mL_{m'}\{ [\lambda +m+m']_q[-\lambda +m+m'+r+s-2]_q\} ^{1/2}
\nonumber\\
 \qquad{}
\times |m+1,k,j,\dots ;m'+1,k',j',\dots \rangle '
\nonumber\\
 \qquad{}
-K_mL_{m'-1}\{ [\lambda +m-m'-s+2]_q[-\lambda +m-m'+r]_q\} ^{1/2}
\nonumber\\  \qquad{}
\times |m+1,k,j,\dots ;m'-1,k',j',\dots \rangle '
\nonumber\\   \qquad{}
+K_{m-1}L_{m'}\{ [\lambda -m+m'-r+2]_q[-\lambda -m+m'+s]_q\} ^{1/2}
\nonumber\\
  \qquad{}\times |m-1,k,j,\dots ;m'+1,k',j',\dots \rangle '
\nonumber\\    \qquad{}-K_{m-1}L_{m'-1}\{ [\lambda -m-m'-r-s+4]_q[-\lambda -m-m'+2]_q\} ^{1/2}
\nonumber\\
  \qquad{}
\times |m-1,k,j,\dots ;m'-1,k',j',\dots \rangle '.\label{24}
\end{gather}

Formulas \eqref{18} and \eqref{24} are used to select $*$-representations in the set
of all irreducible representations $T_{\epsilon \lambda }$ of the algebra
${\rm so}'_q(r,s)$. This selection is fulf\/illed in the same way as in the
case of the $q$-deformed algebras ${\rm so}'_q(2,1)$ and
${\rm so}'_q(3,1)$ in \cite{14}, that is, by a direct check that the relations (8)
are satisf\/ied. This selection leads to the following theorem.

\begin{theorem} All irreducible representations $T_{\epsilon \lambda }$
of ${\rm so}'_q(r,s)$ with $\lambda =-{\overline \lambda} +r+s-2$ are
$*$-rep\-resentations (the principal degenerate series of $*$-representations).
All irreducible representations~$T_{\epsilon \lambda }$ with ${\rm Im}\, \lambda
=\pi /h$  are $*$-representations
(the strange series). If $r$ and $s$ are both even or both odd,
then all irreducible representations
$T_{0\lambda }$, ${1\over 2}(r+s)-1<\lambda
<{1\over 2}(r+s)$, for even ${1\over 2}(r+s)$ and
all irreducible representations
$T_{1\lambda }$, ${1\over 2}(r+s)-1<\lambda
<{1\over 2}(r+s)$, for odd ${1\over 2}(r+s)$ are $*$-representations (the
supplementary series). If among the integers $r$ and $s$ one is odd and
another is even, then all irreducible representations $T_{\epsilon \lambda }$,
${1\over 2}(r+s)-1<\lambda <{1\over 2}(r+s-1)$ are $*$-representations
(supplementary series).
\end{theorem}

   This theorem describes all $*$-representations in the set
of irreducible representations $T_{\epsilon \lambda }$. However, there are
equivalent representations in the formulation of Theorem 2. All possible
equivalences are given by equivalence relations described above
or by their combinations (pro\-ducts).

\section[Reducible representations $T_{\epsilon \lambda}$]{Reducible representations $\boldsymbol{T_{\epsilon \lambda}}$}

 Let us study a structure of reducible representations
$T_{\epsilon \lambda }$ of the algebra ${\rm so}'_q(r,s)$.
Vanishing of some coef\/f\/icients in formula \eqref{18} or \eqref{24} leads to
appearing of invariant subspaces in the carrier space of the representation
$T_{\epsilon \lambda }$. Analysis of reducibility and f\/inding of all
irreducible constituents in $T_{\epsilon \lambda }$ are done in the
same way as in the classical case \cite{11,12,13}. For this reason,
we shall formulate the results of such analysis without detailed proof.

   Let us also note that, as in the classical case, reducible representations
$T_{\epsilon \lambda }$ and $T_{\epsilon ,-\lambda +r+s-2}$ contain the same
(equivalent) irreducible constituents. This is easily seen from formula \eqref{24}.
Studying the representations $T_{\epsilon \lambda }$, we have to
distinguish the cases of odd and even $r$ and $s$ since in dif\/ferent cases
a structure of reducible representations $T_{\epsilon \lambda }$
is dif\/ferent. Below, we investigate all reducible representations
$T_{\epsilon \lambda }$ (which are excluded in Theorem 1).

\subsection[The case of even $r$ and $s$]{The case of even $\boldsymbol{r}$ and $\boldsymbol{s}$}

Let $\lambda$ be an even integer in $T_{0\lambda }$ and
an odd integer in $T_{1\lambda }$.
If $\lambda \le 0$ then in the
carrier space ${\mathcal{H}}(T_{\epsilon \lambda })$ of the representation
$T_{\epsilon \lambda }$ there exist invariant subspaces
\begin{gather*}
{\mathcal{H}}^F_{\lambda }=\bigoplus _{m+m'\le -\lambda }
{\mathcal{V}}(m,{\bf 0} ;m',{\bf 0}),
\qquad
{\mathcal{H}}^0_{\lambda }=\bigoplus _{\lambda -r+2\le m-m'\le -\lambda +s-2}
 {\mathcal{V}}(m,{\bf 0} ;m',{\bf 0}),
\end{gather*}
where ${\mathcal{V}}(m,{\bf 0} ;m',{\bf 0}),$ is the subspace
of ${\mathcal{H}}(T_{\epsilon \lambda })$, on which the irreducible
representation of the subalgebra ${\rm so}'_q(r)\times {\rm so}'_q(s)$
with highest weight $(m,{\bf 0} ;m',{\bf 0})$ is realized. On the subspace
${\mathcal{H}}^F_{\lambda }$ the f\/inite dimensional irreducible representation of
the algebra ${\rm so}'_q(r,s)$ with highest weight $(-\lambda , {\bf 0})$
is realized. An irreducible representation of ${\rm so}'_q(r,s)$
is realized on the quotient space
${\mathcal{H}}^0_{\lambda}/{\mathcal{H}}^F_{\lambda }$. We denote
it by $T^0_\lambda$. A direct sum of two irreducible representations of
${\rm so}'_q(r,s)$ is realized on the quotient space
${\mathcal{H}}(T_{\epsilon \lambda }) /
{\mathcal{H}}^0_\lambda$. One of them acts on the direct sum of the subspaces
${\mathcal{V}}(m,{\bf 0} ;m',{\bf 0})$ for which
$m-m'>-\lambda -s-2$ (we denote it
by $T^-_\lambda $). The second one acts on the direct sum of the subspaces
${\mathcal{V}}(m,{\bf 0} ;m',{\bf 0})$ for which $m'-m>-\lambda +r-2$ (we denote it by
$T^+_\lambda$). Figures showing distribution of subspaces
$\mathcal{V}(m,{\bf 0} ;m',{\bf 0})$ between the subspaces $\mathcal{H}^F_{\lambda }$,
$\mathcal{H}^0_{\lambda }$, $\mathcal{H}^+_{\lambda }$,
$\mathcal{H}^-_{\lambda }$ (and also for other
cases, considered below) are the same as in the
classical case and can be found in~\cite{11}, Chapter 8.

Let now $\lambda$ be even in the representation $T_{0\lambda}$ and odd
in the representation $T_{1\lambda}$, and let $0<\lambda \le
{1\over 2}(r+s)-2$. Then on the space ${\mathcal{H}}(T_{\epsilon \lambda})$ there exists
only one invariant subspace ${\mathcal{H}}^0_\lambda $, which is the orthogonal sum of
the subspaces ${\mathcal{V}}(m,{\bf 0}\ ;m',{\bf 0})$ for which
\[
m-m'\le -\lambda +s-2,\qquad m'-m\le -\lambda +r-2.
\]
The representation of ${\rm so}'_q(r,s)$ on this subspace is irreducible
(we denote it by $T^0_\lambda$). The direct sum of two irreducible
representations of ${\rm so}'_q(r,s)$ is realized on the quotient space
${\mathcal{H}}(T_{\epsilon \lambda })
/{\mathcal{H}}^0_\lambda$. For one of them we have $m-m'>-\lambda +s-2$
(we denote this irreducible representation by~$T^-_\lambda$), and for another
one $m'-m>-\lambda +r-2$ (we denote it by~$T^+_\lambda $). For
$\lambda ={1\over 2}(r+s)-2$ the range of values of $m$ and $m'$ lies
on one line in the coordinate space $(m,m')$. Physicists call such
subrepresentations {\it ladder representations}.

For $\lambda ={1\over 2}(r+s)-1$, the representation $T_{0\lambda}$ of
${\rm so}'_q(r,s)$, if this number $\lambda$ is even, and the
representation $T_{1\lambda}$, if this number is odd,
is a direct sum of two irreducible representations $T^-_\lambda$ and
$T^+_\lambda$: for the f\/irst one we have $m'-m\le -\lambda +r-2$ and for the
second one $m-m'\le -\lambda +s-2$.

   Since the reducible representations $T_{\epsilon \lambda}$
and $T_{\epsilon ,-\lambda +r+s-2}$ contain the same irreducible
constituents, a structure of other reducible
representations in this case is determined by
that of the representations considered above.

   In the same way as in the classical case, it is easy to verify that
the following irreducible representations, considered here, are
$*$-representations:
\smallskip

   (a) all the representations $T^+_\lambda$ and $T^-_\lambda$ (the discrete series);
\smallskip

   (b) the representation $T^0_{(r+s-4)/2}$.
\smallskip

\begin{proposition} Irreducible representations
$T_{\epsilon \lambda }$ and the irreducible representations
$T^+_\lambda$, $T^-_\lambda$, $T^0_\lambda$
of this subsection exhaust all infinite
dimensional irreducible representations
of the algebra ${\rm so}'_q(r,s)$ with even $r$ and $s$
which consist under restriction to ${\rm so}'_q(r)\times
{\rm so}'_q(s)$ of irreducible representations of this
subalgebra only with highest weights of the form
$(m,0,\dots,0)(m',0,\dots,0)$. Moreover, these representations
are pairwise non-equivalent.
\end{proposition}

This proposition is proved in the same way as in the case
of the group $SO_0(r,s)$ (see Chapter~8 in \cite{12}).

\subsection[The case of even $r$ and odd $s$]{The case of even $\boldsymbol{r}$ and odd $\boldsymbol{s}$}

Let $\lambda$ be a non-positive
integer of the same evenness as $m+m'$ does. Then in
the space ${\mathcal{H}}(T_{\epsilon \lambda })$ of the reducible representation
$T_{\epsilon \lambda }$ there exist two invariant subspaces
\begin{gather*}
{\mathcal{H}}^F_{\lambda }=\bigoplus _{m+m'\le -\lambda }
{\mathcal{V}}(m,{\bf 0} ;m',{\bf 0}),
\qquad
{\mathcal{H}}^0_{\lambda }=\bigoplus _{m'-m\le -\lambda +r-2}
 {\mathcal{V}}(m,{\bf 0} ;m',{\bf 0}).
\end{gather*}
The irreducible f\/inite dimensional representation $T^F_\lambda$ of
${\rm so}_q(r,s)$ with the highest weight $(-\lambda ,{\bf 0})$
is realized in the f\/irst subspace. In the quotient spaces
${\mathcal{H}}^0_\lambda /{\mathcal{H}}^F_\lambda$ and
${\mathcal{H}}(T_{\epsilon \lambda }) /{\mathcal{H}}^0_\lambda$ the irreducible
representations of ${\rm so}'_q(r,s)$ are realized which will be
denoted by $T^1_\lambda$ and $T^+_\lambda$ respectively. So, in this case
the representation $T_{\epsilon \lambda }$ consists of three irreducible
constituents.

If $0<\lambda < {1\over 2}(r+s)-2$ and, besides, $\lambda$ is an integer
of the same evenness as $m+m'$ does, then in ${\mathcal{H}}(T_{\epsilon \lambda })$
there exists only one invariant subspace ${\mathcal{H}}^0_\lambda$ which is a
direct sum of the subspaces ${\mathcal{V}}(m,{\bf 0} ;m',{\bf 0})$ for which
$m'-m\le -\lambda +r-2$. The irreducible representations of the algebra
${\rm so}_q(r,s)$ are realized on ${\mathcal{H}}^0_\lambda$ and
${\mathcal{H}}(T_{\epsilon \lambda })/{\mathcal{H}}^0_\lambda $.
We denote them by $T^1_\lambda$ and $T^+_\lambda$ respectively.

If $\lambda <{1\over 2}(r+s)-2$ and, besides, $\lambda$ is an integer
such that $\lambda \equiv (m+m'+1)$ (mod 2), then only one invariant
subspace exists in ${\mathcal{H}}(T_{\epsilon \lambda })$. This subspace is
\[
{\mathcal{H}}^0_\lambda =\bigoplus _{m-m'\le -\lambda +s-2}
{\mathcal{V}}(m,{\bf 0} ;m',{\bf 0}).
\]
The irreducible representations of the algebra ${\rm so}'_q(r,s)$
are realized on ${\mathcal{H}}^0_\lambda$ and
${\mathcal{H}}(T_{\epsilon \lambda })/\mathcal{H}^0_\lambda $. We denote them by
$T^2_\lambda$ and $T^-_\lambda$, respectively.

   Since the reducible representations $T_{\epsilon \lambda}$
and $T_{\epsilon ,-\lambda +r+s-2}$ contain the same irreducible
constituents, a structure of other reducible
representations in this case is determined by
that of the representations considered above.

In the set of irreducible representations, considered here, only the
representations $T^+_\lambda$ and~$T^+_\lambda$ are $*$-representations.

The case of odd $r$ and even $s$ is considered absolutely in the same way,
interchanging the roles of $r$ and $s$ as well as of $m$ and $m'$. For this
reason, we omit consideration of this case.

\begin{proposition} Irreducible representations
$T_{\epsilon \lambda }$ and the irreducible representations
$T^+_\lambda$, $T^-_\lambda$, $T^1_\lambda$, $T^2_\lambda$
of this subsection exhaust all infinite
dimensional irreducible representations
of the algebra ${\rm so}'_q(r,s)$ with even $r$ and odd $s$
which consist under restriction to ${\rm so}'_q(r)\times
{\rm so}'_q(s)$ of irreducible representations of this
subalgebra only with highest weights of the form
$(m,0,\dots,0)(m',0,\dots,0)$. Moreover, these representations
are pairwise non-equivalent.
\end{proposition}

\subsection[The case of odd $r$ and $s$]{The case of odd $\boldsymbol{r}$ and $\boldsymbol{s}$}

If $\lambda $ is a non-positive
integer such that $\lambda \equiv (m+m')$ (mod 2), then in the space
$\mathcal{H}(T_{\epsilon \lambda })$ of the representation $T_{\epsilon \lambda }$
there exists only one invariant subspace
$\mathcal{H}^F_\lambda$ containing all
those subspaces ${\mathcal{V}}(m,{\bf 0} ;m',{\bf 0})$ for which $m+m'\le -\lambda$.
On this invariant subspace the irreducible f\/inite dimensional representation
of ${\rm so}'_q(r,s)$ with the highest weight $(-\lambda ,{\bf 0})$
is realized. On the quotient space ${\mathcal{H}}
(T_{\epsilon \lambda })/\mathcal{H}^F_\lambda$
the irreducible representation of
${\rm so}'_q(r,s)$ is realized which is denoted by $T^3_\lambda$.

If $\lambda \le {1\over 2}(r+s)-2$ and, besides, $\lambda$ is an integer
such that $\lambda \equiv (m+m'+1)$ (mod 2), then in the space
${\mathcal{H}}(T_{\epsilon \lambda })$ there exists only one invariant subspace
${\mathcal{H}}^0_\lambda$ containing all the subspaces
${\mathcal{V}}(m,{\bf 0} ;m',{\bf 0})$ for which
\[
m'-m\le -\lambda +r-2,\qquad m-m'\le -\lambda +s-2.
\]
An irreducible representation of ${\rm so}'_q(r,s)$ is realized on
${\mathcal{H}}^0_\lambda$ (we denote it by $T^0_\lambda$). On the quotient space
${\mathcal{H}}(T_{\epsilon \lambda })
/{\mathcal{H}}^0_\lambda$ a direct sum of two irreducible
representations of ${\rm so}'_q(r,s)$ acts. For one of them we have
$m'-m\le -\lambda +r-2$, and for the other $m-m'\le -\lambda +s-2$. We
denote these irreducible representations by $T^-_\lambda$ and $T^+_\lambda$,
respectively.

If $\lambda ={1\over 2}(r+s)-1$, then the representation
$T_{0\lambda }$ for odd ${1\over 2}(r+s)-1$ and the representation
$T_{1\lambda }$ for even ${1\over 2}(r+s)-1$ decompose into a
direct sum of two irreducible representations of ${\rm so}'_q(r,s)$
(we denote them by $T^-_\lambda$ and $T^+_\lambda$). For the f\/irst
representation we have $m'-m\le -\lambda +r-2$ and for the second one
$m-m'\le -\lambda +s-2$.

   Since the reducible representations $T_{\epsilon \lambda}$
and $T_{\epsilon ,-\lambda +r+s-2}$ contain the same irreducible
constituents, a structure of other reducible
representations in this case is determined by
that of the representations considered above.

   In the set of irreducible representations, considered in this subsection,
only the following ones are $*$-representations:
\smallskip

   (a) all the representations $T^+_\lambda$ and $T^-_\lambda$ (the discrete series);
\smallskip

   (b) the representation $T^0_{(r+s-4)/2}$.
\smallskip

\begin{proposition} Irreducible representations
$T_{\epsilon \lambda }$ and the irreducible representations
$T^+_\lambda$, $T^-_\lambda$, $T^0_\lambda$, $T^3_\lambda$
of this subsection exhaust all infinite
dimensional irreducible representations
of the algebra ${\rm so}'_q(r,s)$ with odd $r$ and $s$
which consist under restriction to ${\rm so}'_q(r)\times
{\rm so}'_q(s)$ of irreducible representations of this
subalgebra only with highest weights of the form
$(m,0,\dots,0)(m',0,\dots,0)$. Moreover, these representations
are pairwise non-equivalent.
\end{proposition}

\pdfbookmark[1]{References}{ref}
\LastPageEnding


\begin{thebibliography}{99}

\footnotesize\itemsep=0pt

\bibitem{1}Gavrilik A.M., Klimyk A.U., $q$-deformed orthogonal and
pseudo-orthogonal algebras and their representations, {\it Lett. Math. Phys.}
{\bf 21} (1991), 215--220.

\bibitem{2}Drinfeld V.G., Hopf algebras and quantum Yang--Baxter equation,
{\it Sov. Math. Dokl.} {\bf 32} (1985), 254--258.

\bibitem{3}Jimbo M., A $q$-dif\/ference analogue of $U_q({\rm gl}(N+1))$ and
the Yang--Baxter equations, {\it Lett. Math. Phys.} {\bf 10} (1985), 63--69.

\bibitem{KS}Klimyk A.U., Schm\"udgen K., Quantum groups and their
representations, Springer, Berlin, 1997.


\bibitem{4}Klimyk A.U., Kachurik I.I.,  Spectra,
eigenvectors and overlap functions for representation operators of
$q$-deformed algebras, {\it Comm. Math. Phys.} {\bf 175} (1996), 89--111.

\bibitem{5}Nelson J., Regge T., 2+1 gravity for genus $s>1$,
{\it Comm. Math. Phys.} {\bf 141} (1991), 211--223.

\bibitem{6}Noumi M., Macdonald's symmetric polynomials as zonal
spherical functions on quantum homogeneous spaces, {\it Adv. Math.} {\bf 123}
(1996), 16--77.

\bibitem{7}Noumi M., Umeda T., Wakayama M.,
 Dual pairs, spherical harmonics and a Capelli identity
in quantum group theory, {\it Compos. Math.} {\bf 104} (1996), 227--277.

\bibitem{8}Iorgov N.Z., Klimyk A.U., The $q$-Laplace
operator and $q$-harmonic polynomials on the quantum vector space,
{\it J. Math. Phys.} {\bf 42} (2001), 1326--1345.

\bibitem{9}Bullock D., Przytycki J.H., Multiplicative
structure of Kauf\/fman bracket skein module quantization,
\mbox{\href{http://arxiv.org/abs/q-alg/9902117}{q-alg/9902117}}.

\bibitem{11a}Twietmeyer E., Real forms of $U_q(g)$, {\it
Lett. Math. Phys.} {\bf 49} (1992), 49--58.

\bibitem{12a}Dobrev V.K., Canonical $q$-deformation of noncompact
Lie (super)algebras, {\it J. Phys. A: Math. Gen.} {\bf 26} (1993),
1317--1329.


\bibitem{13a}Celegini E., Giachetti R., Reyman A., Sorace E., Tarlini M.,
$SO_q(n+1,n-1)$ as a real form of $SO_q(2n,\mathbb{C})$,
{\it Lett. Math. Phys.} {\bf 23} (1991), 45--44.


\bibitem{10}Raczka R., Limic N., Niederle J., Discrete degenerate representations
of the noncompact rotation groups,
{\it J. Math. Phys.} {\bf 7} (1966), 1861--1876.

\bibitem{11}Molchanov V.F., Representations of pseudo-orthogonal groups
associated with a cone, {\it Math. USSR Sbornik} {\bf 10} (1970), 353--347.

\bibitem{12}Klimyk A.U., Matrix elements and Clebsch--Gordan
coef\/f\/icients of group representations, Naukova Dumka, Kiev, 1979.

\bibitem{13}Howe R.E., Tan E.C., Homogeneous functions on light cone:
the inf\/initesimal structure of some degenerate principal series
representations, {\it Bull. Amer. Math. Soc.} {\bf 28} (1993), 1--74.

\bibitem{14}Gavrilik A.M., Klimyk A.U., Representations of
$q$-deformed algebras $U_q({\rm so}_{2,1})$ and $U_q({\rm so}_{3,1})$,
{\it J. Math. Phys.} {\bf 35} (1994), 3670--3686.

\bibitem{15}Kachurik I.I., Klimyk A.U., Representations of the
$q$-deformed algebra $U_q({\rm so}_{r,2})$, {\it Dokl. Akad. Nauk
Ukrainy, Ser. A} (1995), no.~9, 18--20.

\bibitem{17}Schm\"udgen K., Unbounded operator algebras and
representation theory, Birkh\"auser, Basel, 1990.

\bibitem{18}Ostrovskyi V., Samoilenko Yu., Introduction to the theory
of representations of f\/initely presented $*$-algebras,
{\it Reviers in Math. and Math. Phys.} {\bf 11} (1999), 1--261.


\end{thebibliography}
\end{document}